\documentclass[12pt]{article}

\usepackage{CJKutf8}
\usepackage{cite}
\usepackage{graphicx}
\usepackage{amsmath}
\usepackage{amsthm}

\usepackage{geometry}
\usepackage{framed}
\usepackage{multirow}
\usepackage{rotating}
\usepackage{amsfonts}
\usepackage{longtable}
\usepackage{array} 

\usepackage{fancyhdr}
\begin{document}

\title{Conjecture About Arbitrary Even-Order Convexity of $w$-Optimization of the Split CIF}
\date{}

\author{Hao Li 
\thanks{Namely \begin{CJK}{UTF8}{gbsn}李颢\end{CJK} .}
}

\maketitle

\begin{abstract}
The $w$-optimization problem is a mathematical problem abstracted from a useful tool for general data fusion rooted in various engineering tasks. Purely from the perspective of practical applications for this useful tool namely the split covariance intersection filter (Split CIF), it is sufficient to know that the $w$-optimization problem enjoys not only the second-order convexity (namely convexity in conventional sense) but also the fourth-order convexity, thanks to which a guaranteed fast implementation of the Split CIF can be realized. On the other hand, based on certain observations and analysis, the author proposes a conjecture concerning the $w$-optimization problem for further mathematical study. The conjecture is that the $w$-optimization problem has arbitrary even-order convexity, or in other words, it has the second-order convexity, the fourth-order convexity, the sixth-order convexity, and so on --- As mentioned above, whether the $w$-optimization problem has higher-order convexity or not might be of little interest itself for practical applications, yet the author believes that the conjecture would stimulate study on systematic mathematical techniques that are potentially interesting and valuable --- This paper presents the conjecture in details.
\end{abstract}

\section{Introduction}

The split covariance intersection filter (or Split CIF for short) \cite{Li2013a} is a useful tool for general data fusion and can be applied in a variety of engineering tasks such as in cooperative intelligent systems \cite{Li2013b, Wanasinghe2014, Pierre2018, ChenX2020, Li2022TITS, Li2024ITS, Li2025VTT} and in single intelligent system operation \cite{Li2013d, Allig2022, Li2022RAL, Li2023AprilTagNavigation, Li2024IV} as well. An indispensable optimization step involved in the Split CIF, namely \textbf{$w$-optimization} (refer to Section \ref{sec:woptprob} for details of the problem statement), concerns the performance and implementation efficiency of the Split CIF. 

The $w$-optimization problem of the Split CIF enjoys the second-order convexity (namely convexity in conventional sense), a proof of which is provided in the author's works \cite{Li2022FARET_2, Li2022FARET_1}. The $w$-optimization problem further enjoys an even more desirable property namely the fourth-order convexity, thanks to which a guaranteed fast implementation of the Split CIF can be realized for practical applications. Relevant details are provided in the author's another paper \cite{Li2026GF}.

Based on certain observations and analysis, the author proposes a conjecture concerning the $w$-optimization problem, namely that the $w$-optimization problem has arbitrary even-order convexity, or in other words, it has the second-order convexity, the fourth-order convexity, the sixth-order convexity, and so on.

Purely from the perspective of practical applications for the Split CIF, it is sufficient to know that the $w$-optimization problem enjoys not only the second-order convexity but also the fourth-order convexity, and whether the $w$-optimization problem has higher-order convexity or not might be of little interest itself. However, the author believes that the conjecture would stimulate study on systematic mathematical techniques that are potentially interesting and valuable. This paper presents the conjecture in details and explain why the conjecture is difficult to handle.

\section{Conjecture concerning the $w$-optimization problem} \label{sec:conjecture_woptprob}

\subsection{The $w$-optimization problem}  \label{sec:woptprob}

Matrices mentioned in this paper are symmetric matrices by default. Given matrices $\mathbf{P}_{1d}$, $\mathbf{P}_{1i}$, $\mathbf{P}_{2d}$, and $\mathbf{P}_{2i}$ that are positive semi-definite, i.e., $\mathbf{P}_{1d} \geq \mathbf{0}$, $\mathbf{P}_{1i} \geq \mathbf{0}$, $\mathbf{P}_{2d} \geq \mathbf{0}$, $\mathbf{P}_{2i} \geq \mathbf{0}$. Besides, the matrices $\mathbf{P}_{1d} + \mathbf{P}_{1i}$ and $\mathbf{P}_{2d} + \mathbf{P}_{2i}$ which normally correspond to covariances of certain estimates are always positive definite, i.e., $\mathbf{P}_{1d} + \mathbf{P}_{1i} > 0$ and $\mathbf{P}_{2d} + \mathbf{P}_{2i} > 0$. For $w \in [0,1]$, define
\begin{align} \label{eq:defineP}
\mathbf{P}_{1}(w) &= \mathbf{P}_{1d}/w + \mathbf{P}_{1i}, \nonumber \\
\mathbf{P}_{2}(w) &= \mathbf{P}_{2d}/(1-w) + \mathbf{P}_{2i}, \nonumber \\
\mathbf{P}(w) &= (\mathbf{P}_{1}(w)^{-1} + \mathbf{P}_{2}(w)^{-1})^{-1}.
\end{align}
When $w=0$ or $w=1$, $\mathbf{P}(w)$ denotes the limit value as $w \to 0$ or $w \to 1$ respectively. Then we have
\begin{subequations}  \label{eq:P_pos_def}
\begin{align}
\mathbf{P}_{1}(w) &\geq \mathbf{P}_{1d} + \mathbf{P}_{1i} > 0, \\
\mathbf{P}_{2}(w) &\geq \mathbf{P}_{2d} + \mathbf{P}_{2i} > 0, \\
\mathbf{P}(w) &> 0,
\end{align}
\end{subequations}
where the third equation (\ref{eq:P_pos_def}c) is a natural consequence of (\ref{eq:P_pos_def}a) and (\ref{eq:P_pos_def}b).

The $w$-optimization problem involved in the Split CIF \cite{Li2013a} can be formalized as
\begin{equation}  \label{eq:wopt_org}
w = \arg \min_{w \in [0,1]} \det(\mathbf{P}(w))
\end{equation}
or somehow more desirably as
\begin{equation}  \label{eq:wopt}
w = \arg \min_{w \in [0,1]} \ln \det(\mathbf{P}(w)).
\end{equation}
The $w$-optimization problem can also be formalized as
\begin{equation}  \label{eq:wopt_tr}
w = \arg \min_{w \in [0,1]} tr \{ \mathbf{P}(w) \},
\end{equation}
where the trace function $tr \{ \cdot \}$ is used instead of the determinant function $\det(\cdot)$ --- It is worth noting that the determinant and trace functions might be replaced by a generic positive definite real-value matrix function, only if it serves as an appropriate metric (or quasi-metric) that characterizes the ``size'' of the matrix $\mathbf{P}(w)$, yet the two kinds of functions are the popular choices in practical applications
\footnote{The author's another paper \cite{Li2026GF} explains why the determinant function $\det(\cdot)$ is even preferred over the trace function $tr \{ \cdot \}$ for the $w$-optimization problem in practical applications, yet for mathematical study, it is still worth taking (\ref{eq:wopt_tr}) into account as well.}.

\subsection{The arbitrary even-order convexity conjecture}  \label{sec:conjecture_AEOC}

For a generic function $f(x)$, if its $m$-th-order derivative is always non-negative (or positive semi-definite), namely
\begin{equation}  \label{eq:m_order_convexity}
\frac{d^m}{dx^m} f(x) \geq 0,
\end{equation}
then it is said to have the \textbf{$m$-th-order convexity} \cite{Li2026GF}. If $f(x)$ serves as the objective function of an optimization problem, the author also abuses the expression to say that the optimization problem has the $m$-th-order convexity. According to (\ref{eq:m_order_convexity}), the conventional-sense convexity is right the second-order convexity.

The proposed conjecture is that \textbf{the $w$-optimization problem has arbitrary even-order convexity}, more specifically, for $k \in \{1, 2, 3, \cdots\}$ we always have
\begin{subequations}  \label{eq:conjecture_AEOC}
\begin{align}
\frac{d^{2 k}}{d w^{2 k}} \ln \det(\mathbf{P}(w)) &\geq 0,  \\
\frac{d^{2 k}}{d w^{2 k}} tr \{ \mathbf{P}(w) \} &\geq 0,
\end{align}
\end{subequations}
where (\ref{eq:conjecture_AEOC}a) and (\ref{eq:conjecture_AEOC}b) correspond to (\ref{eq:wopt}) and (\ref{eq:wopt_tr}) respectively.

\section{Observations and analysis}

\subsection{Proof of the conjecture for special cases}  \label{sec:conjecture_special_case}

The proposed conjecture stems from the $w$-optimization problem involved in the Split CIF. As explained in \cite{Li2013a}, the Split CIF has two special cases naturally. Here is the first special case: If both matrices $\mathbf{P}_{1d}$ and $\mathbf{P}_{2d}$ are set to zero, then the Split CIF is reduced to the famous Kalman filter \cite{Kalman1960}. As the assumption that $\mathbf{P}_{1d} + \mathbf{P}_{1i} > 0$ and $\mathbf{P}_{2d} + \mathbf{P}_{2i} > 0$ is always expected, the matrices $\mathbf{P}_{1i}$ and $\mathbf{P}_{2i}$ in such special case should be positive definite. Then we have a reduced version of (\ref{eq:defineP}), i.e.
\begin{align} \label{eq:defineP_Pd=0_KF}
\mathbf{P}_{1}(w)^{-1} &= \mathbf{P}_{1i}^{-1}, \nonumber \\
\mathbf{P}_{2}(w)^{-1} &= \mathbf{P}_{2i}^{-1}, \nonumber \\
\mathbf{P}(w)^{-1} &= \mathbf{P}_{1i}^{-1} + \mathbf{P}_{2i}^{-1} \iff \mathbf{P}(w) = (\mathbf{P}_{1i}^{-1} + \mathbf{P}_{2i}^{-1})^{-1},
\end{align}
where the optimization variable $w$ does not take any effect and the matrix $\mathbf{P}(w)$ can be regarded as a constant function in terms of $w$. So
\begin{subequations}  \label{eq:P_1stSpecial_dif_general}
\begin{align}
\frac{d^m}{d w^m} \ln \det(\mathbf{P}) &= 0,  \\ 
\frac{d^m}{d w^m} tr \{ \mathbf{P} \} &= 0,
\end{align}
\end{subequations}
and the conjecture (\ref{eq:conjecture_AEOC}) holds trivially for the first special case.

Here is the second special case: If both matrices $\mathbf{P}_{1i}$ and $\mathbf{P}_{2i}$ are set to zero, then the Split CIF is reduced to the covariance intersection (filter) \cite{Julier1997b}. As the assumption that $\mathbf{P}_{1d} + \mathbf{P}_{1i} > 0$ and $\mathbf{P}_{2d} + \mathbf{P}_{2i} > 0$ is always expected, the matrices $\mathbf{P}_{1d}$ and $\mathbf{P}_{2d}$ in such special case should be positive definite. Then we have a reduced version of (\ref{eq:defineP}), i.e.
\begin{align} \label{eq:defineP_Pi=0}
\mathbf{P}_{1}(w)^{-1} &= \mathbf{P}_{1d}^{-1} w, \nonumber \\
\mathbf{P}_{2}(w)^{-1} &= \mathbf{P}_{2d}^{-1} (1-w), \nonumber \\
\mathbf{P}(w)^{-1} &= \mathbf{P}_{1d}^{-1} w + \mathbf{P}_{2d}^{-1} (1-w) \equiv \mathbf{A} w + \mathbf{B} \iff \mathbf{P}(w) = (\mathbf{A} w + \mathbf{B})^{-1},
\end{align}
where
\begin{align*}
\mathbf{A} \equiv \mathbf{P}_{1d}^{-1} - \mathbf{P}_{2d}^{-1}, \qquad \mathbf{B} \equiv \mathbf{P}_{2d}^{-1} > 0.
\end{align*}

For denotation conciseness in following derivations, explicit writing of ``$(w)$'' for $w$-parametrized variables is omitted. Besides, it is taken for granted that all relevant $w$-parametrized variables are arbitrary-order differentiable. Furthermore, the following equalities
\begin{subequations}  \label{eq:GF_lemmas}
\begin{align}
\frac{d}{dw} \ln \det(\mathbf{M}) &= tr\{\mathbf{M}^{-1} \frac{d \mathbf{M}}{dw} \}  \\
\frac{d \mathbf{M}^{-1}}{dw} &= -\mathbf{M}^{-1} \frac{d \mathbf{M}}{dw} \mathbf{M}^{-1}  \\
tr\{\mathbf{M}_1 \mathbf{M}_2 \cdots \mathbf{M}_k\} &= tr\{\mathbf{M}_2 \mathbf{M}_3 \cdots \mathbf{M}_k \mathbf{M}_1\} = \cdots = tr\{\mathbf{M}_k \mathbf{M}_1 \cdots \mathbf{M}_{k-2} \mathbf{M}_{k-1}\}
\end{align}
\end{subequations}
which correspond to three lemmas used in \cite{Li2026GF} are also used below. Refer to \cite{Li2026GF} for details of the lemmas.

\subsubsection*{Proof of (\ref{eq:conjecture_AEOC}a) for the second special case}

Substitute (\ref{eq:defineP_Pi=0}) into the first two equations of (\ref{eq:GF_lemmas}) and obtain
\begin{subequations}  \label{eq:P_dif_features}
\begin{align}
\frac{d}{dw} \ln \det(\mathbf{P}) &= - \frac{d}{dw} \ln \det(\mathbf{A} w + \mathbf{B}) = - tr\{(\mathbf{A} w + \mathbf{B})^{-1} \frac{d (\mathbf{A} w + \mathbf{B})}{dw} \} = - tr\{ \mathbf{P} \mathbf{A} \},  \\
\frac{d}{dw} \mathbf{P} &= \frac{d}{dw} (\mathbf{A} w + \mathbf{B})^{-1} = - (\mathbf{A} w + \mathbf{B})^{-1} \frac{d (\mathbf{A} w + \mathbf{B})}{dw} (\mathbf{A} w + \mathbf{B})^{-1} = - \mathbf{P} \mathbf{A} \mathbf{P}.
\end{align}
\end{subequations}
Following (\ref{eq:P_dif_features}) and using (\ref{eq:GF_lemmas}c) namely the \textit{cyclic property} of trace operation when necessary in following derivation, we have
\begin{align}  \label{eq:P_2ndSpecial_dif_general}
\frac{d^2}{d w^2} \ln \det(\mathbf{P}) &= tr\{ (\mathbf{P} \mathbf{A})^2 \}, \nonumber \\
\frac{d^3}{d w^3} \ln \det(\mathbf{P}) &= - 2 \cdot tr\{ (\mathbf{P} \mathbf{A})^3 \}, \nonumber \\
\frac{d^4}{d w^4} \ln \det(\mathbf{P}) &= 6 \cdot tr\{ (\mathbf{P} \mathbf{A})^4 \}, \nonumber \\
\cdots & \qquad \cdots \nonumber \\
\frac{d^m}{d w^m} \ln \det(\mathbf{P}) &= (-1)^m (m-1)! \cdot tr\{ (\mathbf{P} \mathbf{A})^m \}.
\end{align}
It is not difficult to arrive at the general formalism (\ref{eq:P_2ndSpecial_dif_general}) via mathematical induction and details are saved here. Then set $m = 2 k$ and obtain
\begin{align*}
\frac{d^{2 k}}{d w^{2 k}} \ln \det(\mathbf{P}) &= (2 k-1)! \cdot tr\{ (\mathbf{P} \mathbf{A})^{2 k} \} = (2 k-1)! \cdot tr\{ \sqrt{\mathbf{P}^{-1}} (\mathbf{P} \mathbf{A})^{2 k} \sqrt{\mathbf{P}} \}  \\
  &= (2 k-1)! \cdot tr\{ (\sqrt{\mathbf{P}} \mathbf{A} \sqrt{\mathbf{P}})^{2 k} \} \geq 0.
\end{align*}

\subsubsection*{Proof of (\ref{eq:conjecture_AEOC}b) for the second special case}

Following (\ref{eq:P_dif_features}b) and using (\ref{eq:GF_lemmas}c) when necessary in following derivation as well, we have
\begin{align}  \label{eq:P_2ndSpecial_dif_general_tr}
\frac{d}{d w} tr \{ \mathbf{P} \} &= - tr\{ \mathbf{P} \mathbf{A} \mathbf{P} \}, \nonumber \\
\frac{d^2}{d w^2} tr \{ \mathbf{P} \} &= 2 \cdot tr\{ (\mathbf{P} \mathbf{A})^2 \mathbf{P} \}, \nonumber \\
\frac{d^3}{d w^3} tr \{ \mathbf{P} \} &= - 6 \cdot tr\{ (\mathbf{P} \mathbf{A})^3 \mathbf{P} \}, \nonumber \\
\frac{d^4}{d w^4} tr \{ \mathbf{P} \} &= 24 \cdot tr\{ (\mathbf{P} \mathbf{A})^4 \mathbf{P} \}, \nonumber \\
\cdots & \qquad \cdots \nonumber \\
\frac{d^m}{d w^m} tr \{ \mathbf{P} \} &= (-1)^m m! \cdot tr\{ (\mathbf{P} \mathbf{A})^m \mathbf{P} \}.
\end{align}
It is not difficult either to arrive at the general formalism (\ref{eq:P_2ndSpecial_dif_general_tr}) via mathematical induction and details are saved here as well. Then set $m = 2 k$ and obtain
\begin{align*}
\frac{d^{2 k}}{d w^{2 k}} tr \{ \mathbf{P} \} &= (2 k)! \cdot tr\{ (\mathbf{P} \mathbf{A})^{2 k} \mathbf{P} \} = (2 k)! \cdot tr\{ \sqrt{\mathbf{P}} (\sqrt{\mathbf{P}} \mathbf{A} \sqrt{\mathbf{P}})^{2 k} \sqrt{\mathbf{P}} \} \geq 0.
\end{align*}
The reason is that $(\sqrt{\mathbf{P}} \mathbf{A} \sqrt{\mathbf{P}})^{2 k}$ is positive semi-definite and hence so is $\sqrt{\mathbf{P}} (\sqrt{\mathbf{P}} \mathbf{A} \sqrt{\mathbf{P}})^{2 k} \sqrt{\mathbf{P}}$.

\subsection{Difficulty of the conjecture for the general case}

Compared with its reduced version for the special cases presented in Section \ref{sec:conjecture_special_case}, the proposed conjecture becomes essentially different in nature for the general case. For the special cases, the matrix $\mathbf{P}(w)^{-1}$, be it given in (\ref{eq:defineP_Pd=0_KF}) or given in (\ref{eq:defineP_Pi=0}), is simple as it is linear in terms of $w$, thanks to which the explicit and tractable formalisms (\ref{eq:P_1stSpecial_dif_general}), (\ref{eq:P_2ndSpecial_dif_general}), and (\ref{eq:P_2ndSpecial_dif_general_tr}) of the generic derivatives $\frac{d^m}{d w^m} \ln \det(\mathbf{P})$ and $\frac{d^m}{d w^m} tr \{ \mathbf{P} \}$ can be derived and accordingly the conjecture for the special cases can be directly verified.

However, for the general case, the matrix $\mathbf{P}(w)^{-1}$ is much more complicated as it becomes nonlinear in terms of $w$. It is at least very hard, if not completely impossible, to derive explicit and tractable formalisms of the generic derivatives $\frac{d^m}{d w^m} \ln \det(\mathbf{P})$ and $\frac{d^m}{d w^m} tr \{ \mathbf{P} \}$ for the general case. By brute-force expansion, the expression of the derivative $\frac{d^m}{d w^m} \ln \det(\mathbf{P})$ quickly becomes too complicated as $m$ increases, which can be seen in \cite{Li2026GF}. The expression of the derivative $\frac{d^m}{d w^m} tr \{ \mathbf{P} \}$ also quickly becomes too complicated as $m$ increases, which can be seen in Appendix \ref{app:dif_OF_tr}. The author believes that the kind of technique presented by himself in \cite{Li2026GF} to prove the fourth-order convexity of the $w$-optimization problem, which relies on brute-force expansion and \textit{ad hoc} mathematical skills, is by no means a correct direction to handle the proposed conjecture.

In a word, we need cleverer and systematic mathematical techniques to handle the proposed conjecture.

\section{Conclusion}

This paper proposes a conjecture concerning the $w$-optimization problem of the Split CIF, namely that the $w$-optimization problem has arbitrary even-order convexity. The conjecture for two typical special cases can be proved, yet the conjecture for the general case is still waiting for systematic mathematical techniques that are potentially interesting and valuable.

\appendix

\section{Derivatives of the trace based objective function}  \label{app:dif_OF_tr}

\subsection{Preliminary expansion}

Compute the first-order, second-order, third-order, and fourth-order derivatives of the trace based objective function $tr \{ \mathbf{P}(w) \}$ consecutively as
\begin{subequations}  \label{eq:trP_Dif}
\begin{align}
&\frac{d}{d w} tr \{ \mathbf{P} \} = tr\{ - (\mathbf{P}_{1}^{-1} + \mathbf{P}_{2}^{-1})^{-1} \frac{d (\mathbf{P}_{1}^{-1} + \mathbf{P}_{2}^{-1})}{d w} (\mathbf{P}_{1}^{-1} + \mathbf{P}_{2}^{-1})^{-1} \} \nonumber \\ 
& \qquad \qquad = tr\{ \mathbf{P} (\mathbf{P}_1^{-1} \frac{d \mathbf{P}_1}{dw} \mathbf{P}_1^{-1} + \mathbf{P}_2^{-1} \frac{d \mathbf{P}_2}{dw} \mathbf{P}_2^{-1}) \mathbf{P} \} \equiv tr\{ \mathbf{P} \mathbf{C}_1 \mathbf{P} \}, \\
&\frac{d^2}{dw^2} tr \{ \mathbf{P} \} = tr\{ 2 \mathbf{P} \mathbf{C}_1 \mathbf{P} \mathbf{C}_1 \mathbf{P} + \mathbf{P} \mathbf{C}_2 \mathbf{P} \}, \\
&\frac{d^3}{dw^3} tr \{ \mathbf{P} \} = tr\{ 6 \mathbf{P} \mathbf{C}_1 \mathbf{P} \mathbf{C}_1 \mathbf{P} \mathbf{C}_1 \mathbf{P} + 3 \mathbf{P} \mathbf{C}_1 \mathbf{P} \mathbf{C}_2 \mathbf{P} + 3 \mathbf{P} \mathbf{C}_2 \mathbf{P} \mathbf{C}_1 \mathbf{P} + \mathbf{P} \mathbf{C}_3 \mathbf{P} \}, \\
&\frac{d^4}{dw^4} tr \{ \mathbf{P} \} = tr\{ 24 \mathbf{P} \mathbf{C}_1 \mathbf{P} \mathbf{C}_1 \mathbf{P} \mathbf{C}_1 \mathbf{P} \mathbf{C}_1 \mathbf{P} + 12 \mathbf{P} \mathbf{C}_1 \mathbf{P} \mathbf{C}_1 \mathbf{P} \mathbf{C}_2 \mathbf{P} + 12 \mathbf{P} \mathbf{C}_1 \mathbf{P} \mathbf{C}_2 \mathbf{P} \mathbf{C}_1 \mathbf{P} \nonumber \\ 
& \qquad + 12 \mathbf{P} \mathbf{C}_2 \mathbf{P} \mathbf{C}_1 \mathbf{P} \mathbf{C}_1 \mathbf{P} + 6 \mathbf{P} \mathbf{C}_2 \mathbf{P} \mathbf{C}_2 \mathbf{P} + 4 \mathbf{P} \mathbf{C}_1 \mathbf{P} \mathbf{C}_3 \mathbf{P} + 4 \mathbf{P} \mathbf{C}_3 \mathbf{P} \mathbf{C}_1 \mathbf{P} + \mathbf{P} \mathbf{C}_4 \mathbf{P} \},
\end{align}
\end{subequations}
where
\begin{align*}
&\mathbf{C}_1 \equiv \mathbf{P}_1^{-1} \frac{d \mathbf{P}_1}{dw} \mathbf{P}_1^{-1} + \mathbf{P}_2^{-1} \frac{d \mathbf{P}_2}{dw} \mathbf{P}_2^{-1}, \\
&\mathbf{C}_2 \equiv \frac{d \mathbf{C}_1}{d w} = -2 \mathbf{P}_1^{-1} \frac{d \mathbf{P}_1}{dw} \mathbf{P}_1^{-1} \frac{d \mathbf{P}_1}{dw} \mathbf{P}_1^{-1} + \mathbf{P}_1^{-1} \frac{d^2 \mathbf{P}_1}{d w^2} \mathbf{P}_1^{-1}  \\
  &\qquad \qquad \qquad - 2 \mathbf{P}_2^{-1} \frac{d \mathbf{P}_2}{dw} \mathbf{P}_2^{-1} \frac{d \mathbf{P}_2}{dw} \mathbf{P}_2^{-1} + \mathbf{P}_2^{-1} \frac{d^2 \mathbf{P}_2}{d w^2} \mathbf{P}_2^{-1}, \\
&\mathbf{C}_3 \equiv \frac{d \mathbf{C}_2}{d w} = 6 \mathbf{P}_1^{-1} \frac{d \mathbf{P}_1}{dw} \mathbf{P}_1^{-1} \frac{d \mathbf{P}_1}{dw} \mathbf{P}_1^{-1} \frac{d \mathbf{P}_1}{dw} \mathbf{P}_1^{-1} - 3 \mathbf{P}_1^{-1} \frac{d^2 \mathbf{P}_1}{d w^2} \mathbf{P}_1^{-1} \frac{d \mathbf{P}_1}{dw} \mathbf{P}_1^{-1}  \\
  &\qquad \qquad \qquad - 3 \mathbf{P}_1^{-1} \frac{d \mathbf{P}_1}{dw} \mathbf{P}_1^{-1} \frac{d^2 \mathbf{P}_1}{d w^2} \mathbf{P}_1^{-1} + \mathbf{P}_1^{-1} \frac{d^3 \mathbf{P}_1}{d w^3} \mathbf{P}_1^{-1}  \\
  &\qquad \qquad \qquad + 6 \mathbf{P}_2^{-1} \frac{d \mathbf{P}_2}{dw} \mathbf{P}_2^{-1} \frac{d \mathbf{P}_2}{dw} \mathbf{P}_2^{-1} \frac{d \mathbf{P}_2}{dw} \mathbf{P}_2^{-1} - 3 \mathbf{P}_2^{-1} \frac{d^2 \mathbf{P}_2}{d w^2} \mathbf{P}_2^{-1} \frac{d \mathbf{P}_2}{dw} \mathbf{P}_2^{-1}  \\
  &\qquad \qquad \qquad - 3 \mathbf{P}_2^{-1} \frac{d \mathbf{P}_2}{dw} \mathbf{P}_2^{-1} \frac{d^2 \mathbf{P}_2}{d w^2} \mathbf{P}_2^{-1} + \mathbf{P}_2^{-1} \frac{d^3 \mathbf{P}_2}{d w^3} \mathbf{P}_2^{-1},  \\
&\mathbf{C}_4 \equiv \frac{d \mathbf{C}_3}{d w} = -24 \mathbf{P}_1^{-1} \frac{d \mathbf{P}_1}{dw} \mathbf{P}_1^{-1} \frac{d \mathbf{P}_1}{dw} \mathbf{P}_1^{-1} \frac{d \mathbf{P}_1}{dw} \mathbf{P}_1^{-1} \frac{d \mathbf{P}_1}{dw} \mathbf{P}_1^{-1} + 12 \mathbf{P}_1^{-1} \frac{d^2 \mathbf{P}_1}{d w^2} \mathbf{P}_1^{-1} \frac{d \mathbf{P}_1}{dw} \mathbf{P}_1^{-1} \frac{d \mathbf{P}_1}{dw} \mathbf{P}_1^{-1} \\
  &\qquad \qquad \qquad + 12 \mathbf{P}_1^{-1} \frac{d \mathbf{P}_1}{dw} \mathbf{P}_1^{-1} \frac{d^2 \mathbf{P}_1}{d w^2} \mathbf{P}_1^{-1} \frac{d \mathbf{P}_1}{dw} \mathbf{P}_1^{-1} + 12 \mathbf{P}_1^{-1} \frac{d \mathbf{P}_1}{dw} \mathbf{P}_1^{-1} \frac{d \mathbf{P}_1}{dw} \mathbf{P}_1^{-1} \frac{d^2 \mathbf{P}_1}{d w^2} \mathbf{P}_1^{-1}  \\
  &\qquad \qquad \qquad -4 \mathbf{P}_1^{-1} \frac{d^3 \mathbf{P}_1}{d w^3} \mathbf{P}_1^{-1} \frac{d \mathbf{P}_1}{dw} \mathbf{P}_1^{-1} -4 \mathbf{P}_1^{-1} \frac{d \mathbf{P}_1}{dw} \mathbf{P}_1^{-1} \frac{d^3 \mathbf{P}_1}{d w^3} \mathbf{P}_1^{-1}  \\
  &\qquad \qquad \qquad -6 \mathbf{P}_1^{-1} \frac{d^2 \mathbf{P}_1}{d w^2} \mathbf{P}_1^{-1} \frac{d^2 \mathbf{P}_1}{d w^2} \mathbf{P}_1^{-1} +4 \mathbf{P}_1^{-1} \frac{d^4 \mathbf{P}_1}{d w^4} \mathbf{P}_1^{-1}  \\
  &\qquad \qquad \qquad -24 \mathbf{P}_2^{-1} \frac{d \mathbf{P}_2}{dw} \mathbf{P}_2^{-1} \frac{d \mathbf{P}_2}{dw} \mathbf{P}_2^{-1} \frac{d \mathbf{P}_2}{dw} \mathbf{P}_2^{-1} \frac{d \mathbf{P}_2}{dw} \mathbf{P}_2^{-1} + 12 \mathbf{P}_2^{-1} \frac{d^2 \mathbf{P}_2}{d w^2} \mathbf{P}_2^{-1} \frac{d \mathbf{P}_2}{dw} \mathbf{P}_2^{-1} \frac{d \mathbf{P}_2}{dw} \mathbf{P}_2^{-1} \\
  &\qquad \qquad \qquad + 12 \mathbf{P}_2^{-1} \frac{d \mathbf{P}_2}{dw} \mathbf{P}_2^{-1} \frac{d^2 \mathbf{P}_2}{d w^2} \mathbf{P}_2^{-1} \frac{d \mathbf{P}_2}{dw} \mathbf{P}_2^{-1} + 12 \mathbf{P}_2^{-1} \frac{d \mathbf{P}_2}{dw} \mathbf{P}_2^{-1} \frac{d \mathbf{P}_2}{dw} \mathbf{P}_2^{-1} \frac{d^2 \mathbf{P}_2}{d w^2} \mathbf{P}_2^{-1}  \\
  &\qquad \qquad \qquad -4 \mathbf{P}_2^{-1} \frac{d^3 \mathbf{P}_2}{d w^3} \mathbf{P}_2^{-1} \frac{d \mathbf{P}_2}{dw} \mathbf{P}_2^{-1} -4 \mathbf{P}_2^{-1} \frac{d \mathbf{P}_2}{dw} \mathbf{P}_2^{-1} \frac{d^3 \mathbf{P}_2}{d w^3} \mathbf{P}_2^{-1}  \\
  &\qquad \qquad \qquad -6 \mathbf{P}_2^{-1} \frac{d^2 \mathbf{P}_2}{d w^2} \mathbf{P}_2^{-1} \frac{d^2 \mathbf{P}_2}{d w^2} \mathbf{P}_2^{-1} +4 \mathbf{P}_2^{-1} \frac{d^4 \mathbf{P}_2}{d w^4} \mathbf{P}_2^{-1}.
\end{align*}

\subsection{Variable transformation}

Define the following notations
\begin{equation}  \label{eq:var_simple}
\mathbf{D}_1 (w) \equiv \frac{\mathbf{P}_{1d}}{w}, \qquad \mathbf{D}_2 (w) \equiv \frac{\mathbf{P}_{2d}}{1-w}, \qquad w_1 \equiv -w, \qquad w_2 \equiv 1-w.
\end{equation}
From definitions given in (\ref{eq:defineP}) we have

\begin{subequations}  \label{eq:P1+P2_difs}
\begin{align}
&\frac{d \mathbf{P}_1}{dw} = \frac{\mathbf{D}_1}{w_1}, \quad \frac{d^2 \mathbf{P}_1}{d w^2} = \frac{2 \mathbf{D}_1}{w_1^2}, \quad \frac{d^3 \mathbf{P}_1}{d w^3} = \frac{6 \mathbf{D}_1}{w_1^3}, \quad \frac{d^4 \mathbf{P}_1}{d w^4} = \frac{24 \mathbf{D}_1}{w_1^4}, \\
&\frac{d \mathbf{P}_2}{dw} = \frac{\mathbf{D}_2}{w_2}, \quad \frac{d^2 \mathbf{P}_2}{d w^2} = \frac{2 \mathbf{D}_2}{w_2^2}, \quad
\frac{d^3 \mathbf{P}_2}{d w^3} = \frac{6 \mathbf{D}_2}{w_2^3}, \quad \frac{d^4 \mathbf{P}_2}{d w^4} = \frac{24 \mathbf{D}_2}{w_2^4}.
\end{align}
\end{subequations}
Substitute (\ref{eq:P1+P2_difs}) into (\ref{eq:trP_Dif}) and obtain
\begin{subequations}  \label{eq:C1+C2+C3+C4}
\begin{align}
&\mathbf{C}_1 = \frac{1}{w_1} \mathbf{P}_1^{-1} \mathbf{D}_1 \mathbf{P}_1^{-1} + \frac{1}{w_2} \mathbf{P}_2^{-1} \mathbf{D}_2 \mathbf{P}_2^{-1},  \\
&\mathbf{C}_2 = \frac{2}{w_1^2} (\mathbf{I} - \mathbf{P}_1^{-1} \mathbf{D}_1) \mathbf{P}_1^{-1} \mathbf{D}_1 \mathbf{P}_1^{-1} + \frac{2}{w_2^2} (\mathbf{I} - \mathbf{P}_2^{-1} \mathbf{D}_2) \mathbf{P}_2^{-1} \mathbf{D}_2 \mathbf{P}_2^{-1},  \\
&\mathbf{C}_3 = \frac{6}{w_1^3} (\mathbf{I} - \mathbf{P}_1^{-1} \mathbf{D}_1)^2 \mathbf{P}_1^{-1} \mathbf{D}_1 \mathbf{P}_1^{-1} + \frac{6}{w_2^3} (\mathbf{I} - \mathbf{P}_2^{-1} \mathbf{D}_2)^2 \mathbf{P}_2^{-1} \mathbf{D}_2 \mathbf{P}_2^{-1},  \\
&\mathbf{C}_4 = \frac{24}{w_1^4} (\mathbf{I} - \mathbf{P}_1^{-1} \mathbf{D}_1)^3 \mathbf{P}_1^{-1} \mathbf{D}_1 \mathbf{P}_1^{-1} + \frac{24}{w_2^4} (\mathbf{I} - \mathbf{P}_2^{-1} \mathbf{D}_2)^3 \mathbf{P}_2^{-1} \mathbf{D}_2 \mathbf{P}_2^{-1},
\end{align}
\end{subequations}
and
\begin{subequations}  \label{eq:trP_Dif2}
\begin{align}
&\frac{d}{d w} tr \{ \mathbf{P} \} = tr\{ \mathbf{P} (\frac{1}{w_1} \mathbf{P}_1^{-1} \mathbf{D}_1 \mathbf{P}_1^{-1} + \frac{1}{w_2} \mathbf{P}_2^{-1} \mathbf{D}_2 \mathbf{P}_2^{-1}) \mathbf{P} \}, \\
&\frac{d^2}{dw^2} tr \{ \mathbf{P} \} = tr\{ 2 [\mathbf{P} (\frac{1}{w_1} \mathbf{P}_1^{-1} \mathbf{D}_1 \mathbf{P}_1^{-1} + \frac{1}{w_2} \mathbf{P}_2^{-1} \mathbf{D}_2 \mathbf{P}_2^{-1})]^2 \mathbf{P}  \\
  &\qquad \qquad + \mathbf{P} [\frac{2}{w_1^2} (\mathbf{I} - \mathbf{P}_1^{-1} \mathbf{D}_1) \mathbf{P}_1^{-1} \mathbf{D}_1 \mathbf{P}_1^{-1} + \frac{2}{w_2^2} (\mathbf{I} - \mathbf{P}_2^{-1} \mathbf{D}_2) \mathbf{P}_2^{-1} \mathbf{D}_2 \mathbf{P}_2^{-1}] \mathbf{P} \}, \\
&\frac{d^3}{dw^3} tr \{ \mathbf{P} \} = tr\{ 6 (\mathbf{P} \mathbf{C}_1)^3 \mathbf{P} + 3 \mathbf{P} \mathbf{C}_1 \mathbf{P} \mathbf{C}_2 \mathbf{P} + 3 \mathbf{P} \mathbf{C}_2 \mathbf{P} \mathbf{C}_1 \mathbf{P} + \mathbf{P} \mathbf{C}_3 \mathbf{P} \}, \\
&\frac{d^4}{dw^4} tr \{ \mathbf{P} \} = tr\{ 24 (\mathbf{P} \mathbf{C}_1)^4 \mathbf{P} + 12 (\mathbf{P} \mathbf{C}_1)^2 \mathbf{P} \mathbf{C}_2 \mathbf{P} + 12 \mathbf{P} \mathbf{C}_1 \mathbf{P} \mathbf{C}_2 \mathbf{P} \mathbf{C}_1 \mathbf{P} \nonumber \\ 
& \qquad + 12 \mathbf{P} \mathbf{C}_2 (\mathbf{P} \mathbf{C}_1)^2 \mathbf{P} + 6 (\mathbf{P} \mathbf{C}_2)^2 \mathbf{P} + 4 \mathbf{P} \mathbf{C}_1 \mathbf{P} \mathbf{C}_3 \mathbf{P} + 4 \mathbf{P} \mathbf{C}_3 \mathbf{P} \mathbf{C}_1 \mathbf{P} + \mathbf{P} \mathbf{C}_4 \mathbf{P} \},
\end{align}
\end{subequations}
where swollen expressions of $\frac{d^3}{dw^3} tr \{ \mathbf{P} \}$ and $\frac{d^4}{dw^4} tr \{ \mathbf{P} \}$ with $\mathbf{C}_1$, $\mathbf{C}_2$, $\mathbf{C}_3$, and $\mathbf{C}_4$ substituted inside are saved.


\newpage
\addcontentsline{toc}{chapter}{Bibliography}

\fancyhf{} 

\bibliographystyle{unsrt}
\bibliography{LI_Hao_Refs_CAEOCWO}

@article{Li2013a,
author={H. Li and F. Nashashibi and M. Yang},
title={Split covariance intersection filter: Theory and its application to vehicle localization},
journal={IEEE Transactions on Intelligent Transportation Systems},
volume={14},
number={4},
pages={1860-1871},
year={2013}
}

@article{Li2013b,
author={H. Li and F. Nashashibi},
title={Cooperative multi-vehicle localization using split covariance intersection filter},
journal={IEEE Intelligent Transportation Systems Magazine},
volume={5},
number={2},
pages={33-44},
year={2013}
}

@conference{Wanasinghe2014, 
author={T. R. Wanasinghe and G. K. I. Mann and R. G. Gosine},
title={Decentralized Cooperative Localization for Heterogeneous Multi-robot System Using Split Covariance Intersection Filter},
booktitle={Canadian Conference on Computer and Robot Vision},
year={2014},
pages={167-174}
}

@conference{Pierre2018, 
author={C. Pierre and R. Chapuis and R. Aufrère and J. Laneurit and C. Debain},
title={Range-Only Based Cooperative Localization for Mobile Robots},
booktitle={International Conference on Information Fusion},
year={2018},
pages={1933-1939}
}

@conference{ChenX2020, 
author={X. Chen and M. Yang and W. Yuan and H. Li and C. Wang},
title={Split Covariance Intersection Filter based Front-Vehicle Track Estimation for Vehicle Platooning without Communication},
booktitle={IEEE Intelligent Vehicles Symposium},
year={2020},
pages={1510-1515}
}

@article{Li2022TITS,
author={S. Fang and H. Li and M. Yang},
title={LiDAR SLAM Based Multivehicle Cooperative Localization Using Iterated {S}plit {CIF}},
journal={IEEE Transactions on Intelligent Transportation Systems},
volume={23},
number={11},
pages={21137-21147},
year={2022}
}

@article{Li2024ITS,
author={S. Fang and H. Li},
title={Multi-vehicle cooperative simultaneous LiDAR {SLAM} and object tracking in dynamic environments},
journal={IEEE Transactions on Intelligent Transportation Systems},
volume={25},
number={9},
pages={11411-11421},
year={2024}
}

@article{Li2025VTT,
author={H. Li and B. Liu and L. Wang},
title={Vehicle top tag assisted vehicle-road cooperative localization for autonomous public buses},
journal={arXiv},
volume={},
number={},
pages={},
year={2025}
}

@conference{Li2013d, 
author={H. Li and F. Nashashibi and B. Lefaudeux and E. Pollard},
title={Track-to-track fusion using split covariance intersection filter-information matrix filter ({SCIF}-{IMF}) for vehicle surrounding environment perception},
booktitle={IEEE International Conference on Intelligent Transportation Systems},
year={2013},
pages={1430-1435}
}

@article{Allig2022,
author={C. Allig and G. Wanielik},
title={Unequal Dimension Track-to-Track Fusion Approaches Using Covariance Intersection},
journal={IEEE Transactions on Intelligent Transportation Systems},
volume={23},
number={6},
pages={5881-5886},
year={2022}
}

@article{Li2022RAL,
author={S. Fang and H. Li and M. Yang and Z. Wang},
title={Inertial navigation system based vehicle temporal relative localization with split covariance intersection filter},
journal={IEEE Robotics and Automation Letters},
volume={7},
number={2},
pages={5270-5277},
year={2022}
}

@article{Li2023AprilTagNavigation,
  title={Split Covariance Intersection Filter Based Visual Localization With Accurate AprilTag Map For Warehouse Robot Navigation},
  author={S. Fang and Y. Li and H. Li},
  journal={arXiv},
  volume={},
  number={},
  pages={},
  year={2023}
}

@article{Li2024IV,
author={Z. Ying and H. Li},
title={{IMM-SLAMMOT}: tightly-coupled {SLAM} and {IMM}-based multi-object tracking},
journal={IEEE Transactions on Intelligent Vehicles},
volume={9},
number={2},
pages={3964-3974},
year={2024}
}

@book{Li2022FARET_1,
author={\begin{CJK}{UTF8}{gbsn}李颢\end{CJK}},
title={\begin{CJK}{UTF8}{gbsn}迭代估计理论基础与应用（英文版）\end{CJK}},
publisher={\begin{CJK}{UTF8}{gbsn}上海交通大学出版社\end{CJK}},
year={2022}
}

@book{Li2022FARET_2,
author={H. Li},
title={Fundamentals and applications of recursive estimation theory},
publisher={Shanghai Jiao Tong University Press},
year={2022}
}

@article{Li2026GF,
  title={Guaranteed Fast Implementation of the Split Covariance Intersection Filter: Nested {N}ewton Method Thanks to the Fourth-Order Convexity of w-Optimization},
  author={H. Li},
  journal={arXiv},
  volume={},
  number={},
  pages={},
  year={2026}
}

@article{Kalman1960,
author={R. Kalman},
title={A new approach to linear filtering and prediction problem},
journal={ASME Trans, Ser. D, J. Basic Eng.},
volume={82},
number={},
pages={35-45},
year={1960}
}

@conference{Julier1997b, 
author={S. Julier and J. Uhlmann},
title={A non-divergent estimation algorithm in the presence of unknown correlations},
booktitle={Proceedings of American Control Conference},
year={1997},
pages={2369-2373}
}

\fancyhead[LE,RO]{\thepage}
\fancyhead[RE]{\textit{ \nouppercase{\leftmark}} }
\fancyhead[LO]{\textit{ \nouppercase{\rightmark}} }

\end{document}